\documentclass[twoside,11pt]{article}
\title{} \author{} \date{}
\usepackage{latexsym,amssymb,times}%,showkeys}
\input amssym.def
\newtheorem{te}{Theorem}[section]

\newtheorem{fac}[te]{Fact}
\newtheorem{lem}[te]{Lemma}

\newtheorem{ex}[te]{Example}
\newtheorem{cor}[te]{Corrolary}
\def\dok{\noindent{\bf Proof. }}
\def\kdok{\hfill $\Box$ \par \vspace*{2mm} }
\def\a{\alpha}

\def\f{\varphi}
\def\o{\omega}
\def\k{\kappa}
\def\r{\rho}
\def\s{\sigma}
\def\x{\xi}
\def\z{\zeta}
\def\t{\tau}
\def\S{{\mathbb S}}
\def\T{{\mathbb T}}
\def\P{{\mathbb P}}
\def\Q{{\mathbb Q}}
\def\B{{\mathbb B}}
\def\N{{\mathbb N}}
\def\K{{\mathbb K}}
\def\X{{\mathbb X}}
\def\Y{{\mathbb Y}}
\def\Z{{\mathbb Z}}

\def\A{{\mathbb A}}
\def\BG{{\mathbb G}}
\def\H{{\mathbb H}}
\def\G{{\mathcal G}}

\def\CB{{\mathcal B}}

\def\CD{{\mathcal D}}
\def\U{{\mathcal U}}
\def\la{\langle}
\def\ra{\rangle}
\def\rw{\rightarrow}
\def\dom{\mathop{\mathrm{dom}}\nolimits}

\def\Emb{\mathop{\rm Emb}\nolimits}

\def\Aut{\mathop{\rm Aut}\nolimits}
\def\Fin{\mathop{\rm Fin}\nolimits}

\def\sm{\mathop{\rm sm}\nolimits}
\def\sq{\mathop{\rm sq}\nolimits}
\def\Pi{\mathop{\rm Pi}\nolimits}
\def\Rado{\mathop{\mathrm{Rado}}\nolimits}
\def\cf{\mathop{\mathrm{cf}}\nolimits}
\def\ps{\mathbin{\emptyset}}
\begin{document}
\thispagestyle{plain}
\begin{center}
          {\large \bf
                    \uppercase{Isomorphic and Strongly Connected Components}}
\end{center}
\begin{center}
{\bf Milo\v s S.\ Kurili\'c\footnote{Department of Mathematics and Informatics, University of Novi Sad,
                                     Trg Dositeja Obradovi\'ca 4, 21000 Novi Sad, Serbia.
                                     e-mail: milos@dmi.uns.ac.rs}
     }
\end{center}
\begin{abstract}
\noindent
We study the partial orderings of the form $\la \P (\X ), \subset\ra$, where $\X$ is a binary relational structure with the connectivity  components
isomorphic to a strongly connected structure $\Y$ and $\P (\X )$ is the set of (domains of) substructures of $\X$ isomorphic to $\X$.
We show that, for example, for a countable $\X$, the poset $\la \P (\X ), \subset\ra$ is either isomorphic to a finite power
of $\P (\Y )$ or forcing equivalent to a separative atomless $\s$-closed poset and, consistently, to $P(\o )/\!\Fin $. %$(P(\o )/\Fin )^+$.
In particular, this holds for  each ultrahomogeneous structure $\X$ such that $\X$ or $\X ^c$ is a disconnected structure
and in this case $\Y$ can be replaced by an  ultrahomogeneous connected digraph.
\end{abstract}
\section{Introduction}\label{S1}
We consider the partial orderings of the form $\la \P (\X), \subset \ra$, where $\X$ is a relational structure and
$\P(\X)$ the set of the domains of its isomorphic substructures.
A rough classification of countable binary structures related to the properties of their posets of copies is
obtained in \cite{Ktow}, defining two structures to be equivalent if the corresponding posets of copies have isomorphic Boolean completions
or, equivalently, are forcing equivalent.
So, for example, for the structures from column $D$ of Diagram 1 of \cite{Ktow} the corresponding posets are forcing equivalent
to an atomless $\omega _1 $-closed poset and, consistently, to $P(\o )/\!\Fin $. This class of structures includes
all scattered linear orders \cite{Kurscatt} (in particular, all countable ordinals \cite{Kurord}),
all structures with maximally embeddable components \cite{Kmaxemb} (in particular, all countable equivalence relations and all disjoint unions of countable ordinals)
and in this paper we show that it contains a large class of ultrahomogeneous structures.

In Theorem \ref{T4024} of Section 3 we show that the poset of copies of a binary structure with $\kappa$-many isomorphic and strongly connected components is either isomorphic to a finite power of the poset of copies of one component, or forcing equivalent to something like  $P(\k )/[\k ]^{<\k }$ and, for countable structures, consistently, to $P(\o )/\!\Fin $.
The main result of Section 4 is that each ultrahomogeneous binary structure which is not biconnected is determined by an ultrahomogeneous digraph in a simple way and this fact is used in Section 5, where we apply Theorem \ref{T4024} to countable ultrahomogeneous binary structures.

\section{Preliminaries}\label{S2}
The aim of this section is to introduce notation and to give basic definitions and facts concerning relational structures and partial orders
which will be used.

We observe {\it binary structures}, the relational structures of the form $\X =\la X, \r \ra$, where $\r $ is a binary relation on the set $X$.
If $\Y =\la Y, \t \ra$ is a binary structure too, a mapping
$f:X \rightarrow Y$ is an {\it embedding} (we write
$f: \X \hookrightarrow  \Y$) iff $f$ is an injection and
$x_1 \r x_2 \Leftrightarrow f(x_1) \t f(x_2)$, for each $x_1, x_2 \in X$.
$\Emb (\X , \Y )$ will denote the set of all embeddings of $\X$ into $\Y$ and, in particular,
$\Emb (\X )  =\Emb (\X , \X ) $.
If, in addition, $f$ is a surjection, $f$ is an {\it  isomorphism} and the structures $\X$ and $\Y$ are called {\it  isomorphic},
in notation $\X \cong \Y$.
If, in particular, $ \Y = \X $, then $f$ is called an {\it  automorphism} of the structure $\X$ and
$\Aut (\X )$ will denote the set of all automorphisms of $\X$.
If $\X = \la X, \r \ra$ is a binary structure, $A\subset X$ and $\r _A = \r \cap (A\times A)$, then $\la A, \r _A \ra$ is the corresponding {\it  substructure}
of $\X$. By $\P (\X )$ we denote the set of domains of substructures of $\X$ which are isomorphic to $\X$, that is
$$
\P (\X )  =  \{ A\subset X : \la A, \r _A \ra \cong \la X, \r \ra\}
          =  \{ f[X] : f \in \Emb (\X )\} .
$$
More generally, if $\X = \la X, \r \ra$ and $\Y =\la Y, \t \ra$ are binary structures we define
$\P (\X , \Y  )  =  \{ B\subset Y : \la B, \t _B \ra \cong \la X, \r \ra\}
                =  \{ f[X] : f \in \Emb (\X , \Y )\} $.
By $\Pi (\X )$ we denote the set of all finite partial isomorphisms of $\X$. A structure $\X$ is called {\it ultrahomogeneous}
iff for each $\f \in \Pi (\X )$ there is $f\in \Aut (\X )$ such that $\f \subset f$.

If $\X _i=\la X_i, \r _i \ra$, $i\in I$, are binary structures  and $X_i \cap X_j =\emptyset$, for
different $i,j\in I$, then the structure $\bigcup _{i\in I} \X _i =\la \bigcup _{i\in I} X_i , \bigcup _{i\in I} \r _i\ra$ will be called
the {\it disjoint union} of the structures $\X _i$, $i\in I$.

If $\la X,\r \ra$ is a binary structure, then
the transitive closure $\r _{rst}$ of the relation $\r _{rs} =\Delta _X \cup \r \cup \r ^{-1}$ (given by $x \;\r_{rst} \;y$ iff there are $n\in \N$ and
$z_0 =x , z_1, \dots ,z_n =y$ such that $z_i \;\r _{rs} \;z_{i+1}$, for each $i<n$)
is the minimal equivalence relation on $X$ containing $\r$.
For $x\in X$ the corresponding element of the quotient $X/\r _{rst}$ will be denoted by $[x]$ and called the
{\it component} of $\la X,\r \ra$ containing $x$. The structure $\la X,\r \ra$ will be called {\it connected}
iff $|X/\r _{rst}|=1$.
It is easy to check (see Proposition 7.2  of \cite{Ktow}) that
$ \la \bigcup _{x\in X}[x],\bigcup _{x\in X}\r_{[x]} \ra$
is the unique representation of $\la X,\r \ra$ as a disjoint union of connected structures.  Also, if $\r ^c = (X\times X)\setminus \r$, then
at least one of the structures $\la X,\r \ra$ and $\la X,\r ^c \ra$ is connected (Proposition 7.3  of \cite{Ktow}). The
following facts (Lemma 7.4 and Theorem 7.5 of \cite{Ktow}) will be used in the sequel.
\begin{fac}\rm\label{T4005}
Let $\la X, \r \ra$ and $\la Y, \t \ra$ be binary structures and $f:X \rightarrow Y$ an embedding. Then
for each $x \in X$

(a) $f[[x]] \subset [f(x)]$;

(b) $f\,|\;[x]:[x]\rightarrow f[[x]]$ is an isomorphism;

(c) If, in addition, $f$ is an isomorphism, then $f[[x]]= [f(x)]$.
\end{fac}
\begin{fac}\rm\label{T4015}
Let $\k$ be a cardinal, let $\X _\a \!= \!\la X_\a , \r _\a \ra , \a \!<\!\k$, be disjoint connected binary structures and $\X$ their union. Then
$C\in \!\P (\X )$ iff there is a function $f\!:\!\k \rightarrow \k$ and there are embeddings $e_\xi : \X _\xi \hookrightarrow \X _{f(\xi )}$, $\xi<\k$,
such that $C= \bigcup _{\xi <\k } e_\xi [X_\xi]$ and
\begin{equation}\label{EQ4009}
\forall \{\xi,\zeta \} \in [\k ]^2 \;\; \forall x \in X_\xi \;\; \forall y \in X_\zeta \;\;\neg \; e_\xi (x)\; \r _{rs} \;e_\zeta (y).
\end{equation}
\end{fac}
Let ${\mathbb P} =\langle  P , \leq \rangle $ be a pre-order. Then
$p\in P$ is an {\it atom}, in notation  $p\in \mathop{\rm At}\nolimits ({\mathbb P} )$, iff each  $q,r\leq p$ are compatible (there is $s\leq q,r$).
${\mathbb P} $ is called {\it atomless} iff $\mathop{\rm At}\nolimits ({\mathbb P} )=\emptyset$; {\it atomic} iff $\mathop{\rm At}\nolimits ({\mathbb P} )$ is dense in ${\mathbb P}$.
If $\kappa $ is a regular cardinal, ${\mathbb P} $ is called $\kappa ${\it -closed} iff for each
$\gamma <\kappa $ each sequence $\langle  p_\alpha :\alpha <\gamma\rangle $ in $P$, such that $\alpha <\beta \Rightarrow p_{\beta}\leq p_\alpha $,
has a lower bound in $P$.
Two pre-orders ${\mathbb P}$ and ${\mathbb Q}$ are called {\it forcing equivalent} iff they produce the same generic extensions. The following fact is folklore.
\begin{fac}\rm\label{T4043}
(a) The direct product of a family of $\k$-closed pre-orders is $\k$-closed.

(b) If $\kappa ^{<\kappa }=\kappa$, then all atomless separative $\kappa$-closed pre-orders  of size $\kappa$ are
forcing equivalent (for example, to the poset (Coll$(\k ,\k ))^+$, or to $(P(\k )/[\k ]^{<\k })^+$).
\end{fac}
A partial order ${\mathbb P} =\langle  P , \leq \rangle $ is called
{\it separative} iff for each $p,q\in P$ satisfying $p\not\leq q$ there is $r\leq p$ such that $r \perp q$.
The {\it separative modification} of ${\mathbb P}$
is the separative pre-order $\mathop{\rm sm}\nolimits ({\mathbb P} )=\langle  P , \leq ^*\rangle $, where
$p\leq ^* q \Leftrightarrow \forall r\leq p \; \exists s \leq r \; s\leq q $.
The {\it separative quotient} of ${\mathbb P}$
is the separative poset $\mathop{\rm sq}\nolimits  ({\mathbb P} )=\langle P /\!\! =^* , \trianglelefteq \rangle$, where
$p = ^* q \Leftrightarrow p \leq ^* q \land q \leq ^* p\;$ and $\;[p] \trianglelefteq [q] \Leftrightarrow p \leq ^* q $.
\begin{fac}    \rm \label{T4042}
(Folklore) Let ${\mathbb P} , {\mathbb Q} $ and ${\mathbb P} _i$, $i\in I$, be partial orderings. Then

(a) ${\mathbb P}$, $\mathop{\rm sm}\nolimits ({\mathbb P})$ and $\mathop{\rm sq}\nolimits  ({\mathbb P})$ are forcing equivalent forcing notions;

(b) ${\mathbb P}$ is atomless iff $\mathop{\rm sm}\nolimits ({\mathbb P} )$ is atomless iff $\mathop{\rm sq}\nolimits  ({\mathbb P} )$ is atomless;

(c) $\mathop{\rm sm}\nolimits ({\mathbb P} )$ is $\kappa $-closed iff $\mathop{\rm sq}\nolimits  ({\mathbb P} )$ is $\kappa $-closed;

(d) ${\mathbb P} \cong {\mathbb Q}$ implies that $\mathop{\rm sm}\nolimits {\mathbb P} \cong  \mathop{\rm sm}\nolimits {\mathbb Q}$ and $\mathop{\rm sq}\nolimits  {\mathbb P} \cong  \mathop{\rm sq}\nolimits  {\mathbb Q}$;

(e) $\mathop{\rm sm}\nolimits (\prod _{i\in I}{\mathbb P} _i) = \prod _{i\in I}\mathop{\rm sm}\nolimits {\mathbb P} _i$ and
    $\mathop{\rm sq}\nolimits  (\prod _{i\in I}{\mathbb P} _i) \cong \prod _{i\in I}\mathop{\rm sq}\nolimits  {\mathbb P} _i$.
\end{fac}
\section{Isomorphic and strongly connected components}\label{S3}
A relational structure $\X =\la X, \r \ra$  will be called {\it  strongly connected}
iff it is connected and for each $A,B\in \P (\X )$ there are $a\in A$ and $b \in B$ such that
$a\; \r _{rs} \; b$. (The structures satisfying $\P (\X )=\{ X\}$ have the second property, but can be disconnected.)
\begin{ex}\rm\label{EX4009}
Some strongly connected structures are: linear orders, full relations, complete graphs,  etc.
The binary tree $\la {}^{<\o }2 , \subset \ra$ is a connected, but not a strongly connected partial order.
\end{ex}

\begin{te}\rm\label{T4024}
Let $\k$ be a cardinal and $\X \!=\bigcup _{\a < \k} \X _\a$ the union of disjoint, isomorphic and strongly connected binary structures.
Then

(a) $\la \P (\X ), \subset \ra \cong \la \P (\X _0) ,\subset \ra ^\k$ and $\sq \la \P (\X ) , \subset  \ra \cong (\sq \la \P (\X _0) ,\subset \ra )^\k$, if $\k<\o$;

(b) $\sq \la \P (\X ) , \subset  \ra$ is an atomless poset, if $\k \geq \o$;

(c) $\sq \la \P (\X ) , \subset  \ra$ is a $\k ^+$-closed poset, if $\k\geq \o$ is regular;

(d) $\sq \la \P (\X ) , \subset  \ra$ is forcing equivalent to the poset $(P(\k )/[\k ]^{<\k })^+$, if
$\k\geq \o$ is regular and $|\P (\X _0 )|\leq 2^\k =\k ^+$. The same holds for $\la \P (\X ) , \subset  \ra$.
\end{te}
\dok
For $A\in [\k ]^{\k }$ and $g\in \prod _{\a \in A }\P (\X _\a )$ let us define $C_g =\bigcup _{\a \in A }g(\a )$.

\vspace{2mm}

\noindent
{\it Claim 1.}
$\P (\X )=\{ C_g : A\in [\k ]^{\k } \land g\in  \prod _{\a \in A }\P (\X _\a )  \} $.

\vspace{2mm}

\noindent
{\it Proof of Claim 1.}
($\subset$) If $C\in \P (\X )$, then, by Fact \ref{T4015}, there is a function $f: \k \rw \k$ and there are embeddings
$e_\xi  : \X _\xi  \hookrightarrow \X _{f(\xi  )}, \xi  < \k$, such that $C=\bigcup _{\xi  \in \k }e_\xi  [X_\xi  ]$ and that (\ref{EQ4009}) is true.

Suppose that $f(\xi  )=f(\zeta  )$, for some different $\xi  , \zeta \in \k$. By the assumption we have $\X _\xi  \cong \X _\zeta  \cong \X _{f(\xi  )}$, which implies
$\P (\X _\xi  ,\X _{f(\xi  )} )=\P (\X _\zeta  ,\X _{f(\xi  )} )=\P (\X _{f(\xi  )} )$. Thus $e_\xi  [X _\xi  ] ,e_\zeta  [X _\zeta  ]\in \P (\X _{f(\xi  )} )$ and,
since the structure $\X _{f(\xi  )}$ is strongly connected, there are $x \in X_\xi  $ and $y \in X_\zeta $ such that $e_\xi  (x) (\r _{f(\xi  )}) _{rs} \; e_\zeta  (y)$,
which, since $\r _{f(\xi  )} \subset \r$, implies $e_\xi  (x)\; \r _{rs} \; e_\zeta  (y)$, which is impossible by (\ref{EQ4009}). Thus $f$ is an injection and, hence,
$A=f[\k ]\in [\k ]^{\k }$. For $f(\xi  )\in f[\k ]$ let $g(f(\xi  )):= e_\xi  [X_\xi  ]$; then
$g(f(\xi  )) \in \P (\X _{f(\xi  )} )$, for all $\xi  \in \k$, that is
$g(\a ) \in \P (\X _\a )$, for all $\a \in A$ and, hence, $g \in \prod _{\a \in A }\P (\X _\a )$. Also
$C=\bigcup _{\xi  \in \k }g(f(\xi  ))=\bigcup _{\a \in A }g(\a )=C_g$ and we are done.

($\supset$) Let $A\in [\k ]^{\k }$, $g\in  \prod _{\a \in A }\P (\X _\a ) $ and let $f:\k \rw A$ be a bijection.
Then for $\xi \in \k$ we have $g(f(\xi ))\in \P (\X _{f(\xi )})=\P (\X _\xi  ,\X _{f(\xi  )} )$ and, hence there is an embedding
$e_\xi  : \X _\xi  \hookrightarrow \X _{f(\xi  )}$ such that $g(f(\xi ))=e_\xi  [X _\xi  ]$. Thus
$C_g =\bigcup _{\a \in A }g(\a )=\bigcup _{\xi  \in \k }g(f(\xi  ))=\bigcup _{\xi  \in \k }e_\xi  [X _\xi  ]$.
If $\xi \neq \zeta  \in \k$, $x \in X_\xi$ and $y \in X_\zeta$, then, since $f$ is an injection, $X_{f(\xi )}$ and $X_{f(\zeta )}$ are different components of
$\X$ containing $e_\xi  (x) $ and  $e_\zeta  (y)$ respectively. So $\neg  e_\xi  (x) \r _{rs}  e_\zeta  (y)$ and (\ref{EQ4009}) is true.
By Fact \ref{T4015} we have $C_g\in \P (\X )$. Claim 1 is proved.
\kdok

(a) By Claim 1 we have $\P (\X )= \{ \bigcup _{i<\k }C_i : \forall i<\k  \;\; C_i \in \P (\X _i)  \}$.  It is easy to see that
the mapping $F$ defined by $F(\la C_i :i<\k \ra )= \bigcup _{i<\k }C_i$ witnesses that the posets $\prod _{i<\k } \la \P (\X _i) ,\subset \ra $
and $\la \P (\X ), \subset \ra$ are isomorphic. Since isomorphic structures have isomorphic posets of copies we have
$\la \P (\X ), \subset \ra \cong \la \P (\X _0) ,\subset \ra ^\k$
and, by Fact \ref{T4042}(d) and (e), $\sq \la \P (\X ), \subset \ra \cong \sq (\la \P (\X _0) ,\subset \ra^\k)\cong (\sq\la \P (\X _0) ,\subset \ra)^\k$.

(b) Let $\k \geq \o$, $\sm \la \P ( \X ) , \subset \ra = \la \P (\X ) , \leq   \ra$
and $\sm \la \P (\X _\a ), \subset \ra =\la \P (\X _\a ), \leq _\a  \ra$, for $\a <\k$. First we prove

\vspace{2mm}

\noindent
{\it Claim 2.}
For each $f,g\in \bigcup _{A\in [\k ]^{\k }} \prod _{\a \in A }\P (\X _\a )$ we have $C_f \leq C_g $ if and only if
\begin{equation}\label{EQ9025}
|(\dom f \setminus \dom g ) \cup \{ \a \in \dom f \cap \dom g : \neg f(\a )\leq _\a g(\a )  \} | <\k ;
\end{equation}

\vspace{2mm}

\noindent
{\it Proof of Claim 2.}
Let $f,g,h\in \bigcup _{A\in [\k ]^{\k }} \prod _{\a \in A }\P (\X _\a )$. Clearly we have
\begin{equation}\label{EQ9027}
C_f \subset C_g \Leftrightarrow \dom f \subset \dom g \land \forall \a \in \dom f \;\; f(\a )\subset g(\a ).
\end{equation}
Let $\perp$ denote the incompatibility relation in the posets $\la \P (\X ), \subset \ra$ and $\la \P (\X _\a), \subset \ra$, $\a <\k$. First we prove
\begin{equation}\label{EQ9028}
C_h \perp C_g \Leftrightarrow |\{ \a \in \dom h \cap \dom g :  h(\a )\not\perp g(\a ) \} |<\k .
\end{equation}
If the set $A=\{ \a \in \dom h \cap \dom g :  h(\a )\not\perp g(\a ) \} $ is of size $\k$, for each $\a \in A$ we choose
$k(\a )\in \P (\X _\a)$ such that $k (\a )\subset h(\a )\cap g(\a )$. So $k \in \prod _{\a \in A }\P (\X _\a )$ and by (a) we have
$C_k \in \P (\X )$. By (\ref{EQ9027}) we have $C_k\subset C_h \cap C_g$ thus $C_h \not\perp C_g$.
Conversely, if $C_h \not\perp C_g$, then by (a) there is $C_k \in \P (\X )$
such that $C_k\subset C_h \cap C_g$. Now $A:= \dom k\in [\k ]^{\k }$ and by (\ref{EQ9027}) we have $A\subset \dom h \cap \dom g $
and $k (\a )\subset h(\a )\cap g(\a )$, for all $\a \in A$. Thus $|\{ \a \in \dom h \cap \dom g :  h(\a )\not\perp g(\a ) \} |=\k$.

Now suppose that $C_f \leq C_g$. Then for each $C_h\in \P (\X )$ satisfying $C_h \subset C_f$  we have $C_h \not\perp C_g $ so, by (\ref{EQ9028}) we have
\begin{equation}\label{EQ9029}
\forall C_h \in \P (\X )\;\; (C_h \subset C_f \Rightarrow |\{ \a \in \dom h \cap \dom g :  h(\a )\not\perp g(\a ) \} |= \k ).
\end{equation}
Suppose that the set $A:=\dom f \setminus \dom g $ is of size $\k$. Then $h:= f\upharpoonright A \in  \prod _{\a \in A }\P (\X _\a )$, clearly $C_h \subset C_f$ and, by (a), $C_h \in \P (\X )$. Also we have $\dom h \cap \dom g =\emptyset$, which is impossible by (\ref{EQ9029}). Thus
\begin{equation}\label{EQ9030}
| \dom f \setminus \dom g | < \k .
\end{equation}
Suppose that the set $A:=\{ \a \in \dom f \cap \dom g : \neg f(\a )\leq _\a g(\a )  \} $ is of size $\k$. For $\a \in A$ there is $C_\a\in \P (\X _\a )$
such that $C_\a \subset f(\a )$ and $C_\a \perp g(\a )$ and we define $h(\a )=C_\a$. Now $h\in \prod _{\a \in A }\P (\X _\a )$,  by (a) we have
$C_h\in \P (\X )$ and, by (\ref{EQ9027}), $C_h \subset C_f$. So by (\ref{EQ9029}) there is $\a \in \dom h \cap \dom g =A$ such that
$C_\a =h(\a )\not\perp g(\a )$, which is not true. Thus
\begin{equation}\label{EQ9031}
|\{ \a \in \dom f \cap \dom g : \neg f(\a )\leq _\a g(\a )  \} | < \k .
\end{equation}
Now from (\ref{EQ9030}) and (\ref{EQ9031}) we obtain (\ref{EQ9025}).

Conversely, assuming (\ref{EQ9030}) and (\ref{EQ9031}) in order to prove $C_f \leq C_g$ we prove (\ref{EQ9029}) first. Let $C_h \in \P (\X )$ and $C_h \subset C_f$. Then, by (\ref{EQ9027}),
\begin{equation}\label{EQ9032}
\dom h \subset \dom f \;\land \;\; \forall \a \in \dom h \;\; h(\a )\subset f(\a ),
\end{equation}
which by (\ref{EQ9030}) implies $| \dom h \setminus \dom g | < \k $ and, hence, $| \dom h \cap \dom g | = \k $.
Since $\dom h \cap \dom g \subset \dom f \cap \dom g$ by (\ref{EQ9031}) we have $|\{ \a \in \dom h \cap \dom g : \neg f(\a )\leq _\a g(\a )  \} | < \k $ and, hence,
$B:= \{ \a \in \dom h \cap \dom g : f(\a )\leq _\a g(\a )  \}$ is a set of size $\k $. By (\ref{EQ9032}), for $\a \in B$ we have $h(\a )\subset f(\a )\leq _\a g(\a )$
which implies $h(\a )\not\perp g(\a )$. So $B \subset \{ \a \in \dom h \cap \dom g :  h(\a )\not\perp g(\a ) \}$ and (\ref{EQ9029}) is true.
Now, by (\ref{EQ9029}) and (\ref{EQ9028}) we have $\forall C_h \in \P (\X )\;\; (C_h \subset C_f \Rightarrow C_h \not\perp C_g )$, that is $C_f \leq C_g$.
Claim 2 is proved.
\kdok

\noindent
Let $A_1$ and $A_2$ be disjoint elements of $[\k ]^\k$. By Claim 1, $C_1=\bigcup _{\a \in A_1}X_\a$ and
$C_2=\bigcup _{\a \in A_2}X_\a$ are disjoint elements of $\P (\X )$ and, hence, they are incompatible in $\la \P (\X ),\subset \ra$.
So, by Theorem 2.2(c) of \cite{Ktow}, the poset $\la \P (\X ),\subset \ra$ is atomless
and, by Fact \ref{T4042}(b), the poset $\sq \la \P (\X ),\subset \ra$ is atomless too.

(c) Let $\k \geq \o$ be a regular cardinal. By Fact \ref{T4042}(c), it is sufficient to prove that the pre-order $\sm \la \P (\X ), \leq\ra$
 is $\k ^+$-closed. Let $\la C_{f_\x } : \x <\k \ra$ be a decreasing sequence in $\la \P (\X ) ,\leq \ra$, that is
\begin{equation}\label{EQ9033}
\forall \z _1 ,\z _2 < \k \;\; (\z _1 <\z _2 \Rightarrow C_{f_{\z _2}} \leq  C_{f_{\z _1}}).
\end{equation}
For $\z _1 ,\z _2 <\k $ let
\begin{equation}\label{EQ9034}
K _{\z _2 ,\z _1 }= \{ \a \in \dom f_{\z _2} \cap \dom f_{\z_1 } : \neg f_{\z _2}(\a )\leq _\a  f_{\z _1 }(\a )  \} .
\end{equation}
Then, by (\ref{EQ9033}) and (c)
\begin{equation}\label{EQ9035}
\forall \z _1 ,\z _2 < \k \;\; (\z _1 <\z _2 \Rightarrow  |\dom f_{\z _2} \setminus  \dom f_{\z_1}|<\k \land |K _{\z _2 ,\z _1 }|<\k )
\end{equation}
and we prove that
\begin{equation}\label{EQ9037}\textstyle
\forall \x<\k \;\; |\bigcap _{\z \leq \x }\dom f_{\z }|= \k .
\end{equation}
First $\bigcap _{\z \leq \x }\dom f_{\z }= \bigcap _{\z < \x }\dom f_{\x } \cap \dom f_{\z }=\dom f_{\x } \cap \bigcap _{\z < \x }(\dom f_{\x }^c \cup \dom f_{\z })$
$= \dom f_{\x } \setminus \bigcup _{\z < \x }(\dom f_{\x } \setminus \dom f_{\z })$.
By (\ref{EQ9035}), $|\dom f_{\x} \setminus  \dom f_{\z}|<\k$, for all $\z <\x$  and, since $|\x |<\k$, by the regularity of $\k$ we have
$|\bigcup _{\z < \x }(\dom f_{\x } \setminus \dom f_{\z })|<\k$ which, since by (a) we have $|\dom f_\x|=\k$, implies (\ref{EQ9037}).

By recursion we define a sequence $\la \a _\x :\x <\k \ra$ in $\k$ as follows.

Let $\a _0 =\min \dom f_0$.

If $\x <\k$ and $\a _\z \in \k$ are defined for $\z<\x$, then for all $\z < \x$ by (\ref{EQ9035}) we have $|K _{\x ,\z  }|<\k $
and, clearly,  $|\a _\z +1 |<\k$ so, by (\ref{EQ9037}) and the regularity of $\k$, we can define
\begin{equation}\label{EQ9036}\textstyle
\a _\x =\min \Big[ \Big( \bigcap _{\z \leq \x }\dom f_\z \Big)\setminus \Big( \bigcup _{\z < \x }K _{\x ,\z  } \cup  \bigcup _{\z < \x }(\a _\z +1 )\Big)\Big] .
\end{equation}
By (\ref{EQ9036}), $\la \a _\x :\x <\k \ra$ is an increasing sequence and, hence, $A:=\{ \a _\x :\x <\k \}\in [\k ]^\k$.
By (\ref{EQ9036}) again, for $\x <\k$ we have $\a _\x \in \dom f_\x$ so  $f_\x (\a _\x )\in \P (\X _{\a _\x })$.
So, for $f\in \prod _{\a _\x \in A}\P (\X _{\a _\x })$, defined by $f (\a _\x )= f_\x (\a _\x )$, for $\x <\k$,
by (a) we have $C_f\in \P (\X )$.

It remains to be shown that for each $\x _0\in \k$ we have $C_f \leq C_{f_{\x _0}}$, that is, by (c),
\begin{equation}\label{EQ9038}
|A \setminus \dom f_{\x _0} | < \k \;\; \mbox{ and}
\end{equation}
\begin{equation}\label{EQ9039}
|\{ \x <\k : \a _\x \in \dom f_{\x _0} \land \neg f_\x (\a  _\x)\leq _{\a _\x} f_{\x _0}(\a  _\x)  \} | < \k .
\end{equation}
By (\ref{EQ9036}), for each $\x\geq \x _0$ we have $\a _\x \in \bigcap _{\z \leq \x }\dom f_\z \subset \dom f_{\x _0}$ and, hence,
$A \setminus \dom f_{\x _0}\subset \{ \a _\x : \x < \x _0 \}$ and (\ref{EQ9038}) is true.

For a proof of (\ref{EQ9039}) it is sufficient to show that
\begin{equation}\label{EQ9040}
\forall \x >\x_0 \;\; f_\x (\a  _\x)\leq _{\a _\x} f_{\x _0}(\a  _\x)  .
\end{equation}
By (\ref{EQ9036}), for $\x >\x_0$ we have $\a _\x \in \dom f_\x \cap \dom f_{\x _0}$ and  $\a _\x \not\in K _{\x ,\x _0}$, that is
$\a _\x \not\in \{ \a \in \dom f_\x \cap \dom f_{\x _0} : \neg f_\x (\a )\leq _\a  f_{\x _0 }(\a )  \}$ thus
$f_\x (\a  _\x)\leq _{\a _\x} f_{\x _0}(\a  _\x)$ and (\ref{EQ9040}) is true.

(d)
Let $\k\geq \o$ be a regular cardinal and $|\P (\X _\a )|\leq 2^\k =\k ^+$, for all $\a <\k$.
Then for $A \in [\k ]^\k$ we have $|\prod _{\a \in A }\P (\X _\a )|\leq (2^\k )^\k =2^\k =\k ^+$ and,
by Claim 1, $|\P (\X )|\leq |\bigcup _{A \in [\k ]^\k }\prod _{\a \in A }\P (\X _\a )|\leq 2^\k 2^\k =2^\k =\k ^+$, which implies
$|\sq \P (\X )|\leq \k ^+$. By (b) and (c) $\sq \P (\X )$ is an atomless $\k ^+$-closed poset and, hence, it contains
a copy of the reversed tree $\la 2^{\leq \k}, \supset \ra$ thus $|\sq \P (\X )|= \k ^+$. (Another way to prove this is to use an
almost disjoint family ${\mathcal A}\subset [\k ]^\k$ of size $\k ^+$; then $\{ \bigcup _{\a\in A}X_\a : A\in {\mathcal A}\}\subset \P (\X )$
determines an antichain in $\sq \P (\X )$ of size $\k ^+$.) Since  $(\k ^+ )^{< \k ^+}=(2^\k)^\k=\k ^+$, by Fact \ref{T4043}(b) the poset
$\sq \P (\X )$ is forcing equivalent to the poset $(P(\k )/[\k ]^{<\k })^+$ (since it is an atomless separative $\kappa ^+$-closed poset of size $\kappa^+$).
By Fact \ref{T4042}(a), the same holds for $\la \P (\X ) , \subset  \ra$.
\hfill $\Box$
\begin{cor}\rm\label{T9022}
If $\k \leq \o $ and $\X =\bigcup _{n < \k} \X _n$ is the union of disjoint, isomorphic and strongly connected binary structures,
then

(a) $\la \P (\X ), \subset \ra \cong \la \P (\X _0) ,\subset \ra ^\k$ and $\sq \la \P (\X ) , \subset  \ra \cong (\sq \la \P (\X _0) ,\subset \ra )^\k$, if $\k<\o$;

(b) If $\k =\o$, then $\sq \la \P (\X ) , \subset  \ra$ is a separative atomless and $\o _1$-closed poset. Under CH it is
    forcing equivalent to the poset $(P(\o )/\Fin )^+$.
\end{cor}
The following examples show that for infinite cardinals $\k$  the statements of Theorem \ref{T4024} are the best possible.

\begin{ex}\rm\label{EX4015}
The posets $\sq \la \P (\X ), \subset \ra$  and $(P(\k )/[\k ]^{<\k })^+$ are not forcing equivalent, although
$\k\geq \o$ is regular and $|\P (\X _\a )|\leq 2^\k $.

Let $\X = \bigcup _{i<\o }\X _i $ be the union of countably many copies $\X _i\!=\!\la X_i, <_i \ra$ of the linear order $\la\o , <\ra$. Then,
since linear orders are strongly connected,
by Theorem \ref{T4024} the poset $\sq \la \P (\X ), \subset \ra$ is atomless, $\o _1$-closed and, clearly, of size $2^\o$.
If, in addition $2^\o =\o _1$, then $\sq\la \P ( \X ) , \subset \ra $ is forcing equivalent to the poset $(P(\o )/\Fin )^+$.

Since, in addition, the components of ${\mathbb X}$ are maximally embeddable (which means that
${\mathbb P} ({\mathbb X} _i , {\mathbb X} _j )=[{\mathbb X} _j ]^{|{\mathbb X} _i|}$, for $i,j\in \o $), by the results of \cite{Kmaxemb}
the poset $\mathop{\rm sq}\nolimits  \langle  {\mathbb P} ({\mathbb X} ), \subset \rangle $ is isomorphic to
the poset $(P(\omega \times \omega)/(\mathop{\rm Fin}\nolimits  \times \mathop{\rm Fin}\nolimits  ))^+$, which is not $\omega _2$-closed \cite{Szym} and,
consistently, neither ${\mathfrak t}$-closed nor ${\mathfrak h}$-distributive \cite{Her}. Thus
in some models of ZFC the posets $\sq \la \P (\X ), \subset \ra$  and $(P(\o )/\Fin )^+$ are not forcing equivalent.
\end{ex}
\begin{ex}\rm\label{EX4016}
In some models of ZFC the poset $\sq \la \P (\X ), \subset \ra$ is not $\k ^{++}$ closed, although the posets
$\sq \la [\k ]^\k , \subset \ra$ and $\sq \la \P (\X _\a ), \subset \ra$, $\a < \k$ are (take $\k =\o$, a model satisfying
${\mathfrak t}>\o _1$ and $\X$ from Example \ref{EX4015}).
\end{ex}
\begin{ex}\rm\label{EX9002}
Statement (c) of Theorem \ref{T4024} is not true for a singular $\k$. It is known that the algebra $P(\k )/[\k ]^{<\k }$ is not
$\o _1$-distributive and, hence, the poset $(P(\k )/[\k ]^{<\k })^+$ is not $\o _2$-closed, whenever $\k$ is a cardinal satisfying $\k > \cf (\k )=\o$ (see \cite{Balc}, p.\ 377). For $\a <\k$ let $\X _\a =\la \{ \a \} , \emptyset\ra$ and let $\X =\bigcup _{\a <\k }\X _\a$. Then it is easy to see that $\P (\X )=[\k ]^\k$ and
$\sq \la \P (\X ), \subset \ra = (P(\k )/[\k ]^{<\k })^+$. Thus the poset $\sq \la \P (\X ), \subset \ra$ is not $\o _2$-closed and, since $\k \geq \aleph _\o$,
it is not $\k ^+$-closed.
\end{ex}
\section{Non biconnected ultrahomogeneous structures}\label{S4}
A binary structure $\X =\la X, \r  \ra $ is a {\it  directed graph (digraph)} iff for each $x,y\in X$ we have
$\neg x\r x$ ($\r$ is irreflexive) and $\neg x\r y \lor \neg y\r x$ ($\r $ is asymmetric). If, in addition,
$x\r y \lor y\r x$, for each different $x,y\in X$, then $\X$ is a {\it  tournament}.
For convenience we introduce the following notation.
If $\X =\la X, \r  \ra $ is a  binary structure, then
its {\it  complement}, $\la X, \r ^c  \ra $, where $\r ^c = X^2 \setminus \r$, will be denoted by $\X ^c$,
its {\it  inverse}, $\la X, \r ^{-1} \ra$, by $\X ^{-1}$,
its {\it  reflexification}, $\la X, \r \cup \Delta _X \ra$, by $\X _{re} $ and
its {\it  irreflexification}, $\la X, \r \setminus \Delta _X \ra$, by $\X _{ir} $.
The binary relation $\r _e$ on $X$ defined by
\begin{equation}\label{EQ9016}
x\r _e y \;\Leftrightarrow \;  x\r y \lor (x\neq y \land \neg x\r y \land \neg y\r x )
\end{equation}
will be called the {\it  enlargement} of $\r$ and the corresponding structure, $\la X, \r _e \ra$, will be denoted by $\X _e$.
A structure $\X$ will be called {\it  biconnected} iff both $\X$ and $\X ^c$ are connected structures.
The following theorem is the main result of this section.
\begin{te}\label{T9017}\rm
For each reflexive or irreflexive ultrahomogeneous binary structure $\X$  we have

- Either $\X$ is biconnected,

- Or there are an  ultrahomogeneous digraph $\Y$ and a cardinal $\k >1$ such that the structure $\X$ is isomorphic to $\bigcup _\k \Y _e$,
$(\bigcup _\k \Y _e)_{re}$, $(\bigcup _\k \Y _e )^c$ or $((\bigcup _\k \Y _e)_{re})^c$.
\end{te}
A proof of Theorem \ref{T9017} is given at the end of this section. It is based on the following statement concerning irreflexive structures.
\begin{te}\label{T9015}\rm
An irreflexive disconnected binary structure is ultrahomogeneous iff its components are isomorphic to the enlargement of an  ultrahomogeneous digraph.
\end{te}
Theorem \ref{T9015} follows from two lemmas given in the sequel.
A binary structure $\X =\la X, \r  \ra $ is called {\it  complete} (see \cite{Fra}, p.\ 393) iff
\begin{equation}\label{EQ9004}
\forall x,y  \;\; (x\neq y \Rightarrow x\r y \lor y\r x ).
\end{equation}
\begin{lem}\label{T9010}\rm
An irreflexive disconnected binary structure $\X$ is ultrahomogeneous iff its components are isomorphic, ultrahomogeneous and complete.
\end{lem}
\dok
Let $\X =\la X , \r \ra= \bigcup _{i\in I}\X _i$, where $\X _i=\la X _i , \r _i\ra$, $i\in I$, are  disjoint, irreflexive and connected binary structures
and $|I|>1$.
%By Fact \ref{T9011}, w.l.o.g.\ we suppose that $\X $ is irreflexive.

($\Rightarrow$)
Suppose that $\X$ is ultrahomogeneous.
Then, for $i,j\in I$, $x\in X_i$ and $y\in X_j$ we have $\f =\{ \la x,y \ra\} \in \Pi (\X )$ and there is $f\in \Aut (\X )$
such that $\f \subset f$. By (c) and (b) of Fact \ref{T4005}, $f|X _i : \X _i \rightarrow \X _j$ is an isomorphism. Thus $\X _i \cong \X _j$.

For $i\in I$ and  $\f  \in \Pi (\X _i)$ we have $\f  \in \Pi (\X )$ and there is $f\in \Aut (\X )$
such that $\f \subset f$. Again, by (c) and (b) of Fact \ref{T4005}, $f|X _i : \X _i \rightarrow \X _i$ is an isomorphism, that is $f|X _i\in \Aut (\X _i)$.
Thus the structure $\X _i$ is ultrahomogeneous.

Suppose that for some $i\in I$ there are
different elements $x$ and $y$ of $X_i$ satisfying $\neg x\r y$ and $\neg y\r x$. Let $j\in I \setminus \{ i \}$ and $z\in X_j$.
Then $\f =\{ \la x,x \ra, \la y,z \ra \} \in \Pi (\X )$ and there is $f\in \Aut (\X )$
such that $\f \subset f$. But then, by Fact \ref{T4005}(c) we would have both $f[X _i] = X _i$ and $f[X _i] = X _j$, which is, clearly, impossible.
Thus the structures $\X _i$ are complete.

($\Leftarrow$) Suppose that the components $\X _i$, $i\in I$, of $\X$ are ultrahomogeneous, isomorphic and complete.
Let $\f \in \Pi (\X)$, where $\dom \f = Y $ and $\f [Y]= Z$, let $J=\{ i\in I : Y\cap X_i \neq\ps \}$ and, for $j\in J$, let $Y_i=Y \cap X_i$ and
$Z_i=\f [Y_i]$. By (\ref{EQ9004}), the structures $\Y _i=\la Y_i , \r _{Y _i}\ra=\la Y_i , (\r _i ) _{Y _i}\ra$, $i\in J$, are connected and, clearly, disjoint,
thus $\Y =\bigcup _{i\in J}\Y _i$ and $\Y _i$, $i\in J$, are the components of $\Y$. Since the restrictions $\f | Y_i : Y_i \rightarrow Z_i$ are isomorphisms,
the structures $\Z _i=\la Z_i , \r _{Z _i}\ra$, $i\in J$, are connected too and, since $\f$ is a bijection, disjoint. Thus
$\Z =\bigcup _{i\in J}\Z _i$ and $\Z _i$, $i\in J$, are the components of $\Z$.

Since $\f : \Y \hookrightarrow \X$, by Fact \ref{T4005}(a) for each $i\in J$ there is $k_i\in I$ such that
$Z_i \subset X_{k_i}$. Suppose that $k_i=k_j=k$, for some different $i,j\in J$. Then, for $x\in Y_i$ and $y\in Y_j$ we would have
$\neg x\r y$ and $\neg y\r x$ and, hence, $\neg \f (x)\r \f (y)$ and $\neg \f (y )\r \f (x)$, which is impossible since $\f (y ),\f (x)\in X_k$
and $\X _k$ satisfies (\ref{EQ9004}).

Thus the mapping $i\mapsto k_i$ is a bijection and there is a bijection $f:I \rightarrow I$ such that $f(i)=k_i$, for all $i\in J$.
Since the structures $\X _i$ are isomorphic, for each $i\in I$ there is an isomorphism $g_i : \X _i \rightarrow \X _{f(i)}$.

For $i\in J$ we have $g_i ^{-1}\circ (\f | Y_i ): \Y _i \hookrightarrow \X _i$ and, hence, $g_i ^{-1}\circ (\f | Y_i )\in \Pi (\X _i)$.
So, since the structure $\X _i$ is ultrahomogeneous, there is $h_i \in \Aut (\X _i)$ such that  $g_i ^{-1}\circ (\f | Y_i )\subset h_i$.
Now $g_i \circ h_i : \X _i \rightarrow \X_{f(i)}$ is an isomorphism and for $x\in Y_i$ we have $g_i (h_i (x))= g_i (g^{-1}(\f (x)))=\f (x)$,
which implies
\begin{equation}\label{EQ9005}
(g_i \circ h_i )|Y_i = \f |Y_i .
\end{equation}
Now it is easy to check that $F=\bigcup _{i\in I\setminus J}g_i \cup \bigcup _{i\in  J}g_i \circ h_i :\X \rightarrow \X$ is an automorphism
of $\X$ and, by (\ref{EQ9005}), $\f \subset F$. Thus $\X$ is an ultrahomogeneous structure.
\kdok
In the sequel we will use the following elementary fact.
\begin{fac}\rm\label{T9011}
Let $\X =\la X, \r  \ra $ be a binary structure. Then

(a) $\Pi (\X )=\Pi ( \X ^c )=\Pi ( \X ^{-1} )$ and $\Aut (\X )=\Aut ( \X ^c )=\Aut ( \X ^{-1} )$;
hence $\X$ is ultrahomogeneous iff $\X ^c$ is ultrahomogeneous iff $\X ^{-1}$ is ultrahomogeneous.
Also $\Emb (\X )=\Emb ( \X ^c )=\Emb ( \X ^{-1} )$; hence $\P (\X )=\P ( \X ^c )=\P ( \X ^{-1} )$.

(b) If $\r $ is an irreflexive relation, then $\Pi (\X )=\Pi ( \X _{re} )$, $\Aut (\X )=\Aut ( \X _{re} )$ and,
hence, $\X$ is ultrahomogeneous iff $\X _{re}$ is ultrahomogeneous.
Also $\Emb (\X )=\Emb ( \X _{re} )$; hence $\P (\X )=\P ( \X _{re} )$.

(c) If $\r $ is a reflexive relation, then $\Pi (\X )=\Pi ( \X _{ir} )$,  $\Aut (\X )=\Aut ( \X _{ir} )$ and,
hence, $\X$ is ultrahomogeneous iff $\X _{ir}$ is ultrahomogeneous.
Also $\Emb (\X )=\Emb ( \X _{ir} )$; hence $\P (\X )=\P ( \X _{ir} )$.

(d) If $\X $ is a digraph, then $\X _e =((\X ^{-1})_{re})^c$. So
$\Pi (\X )=\Pi ( \X _{e} )$, $\Aut (\X )=\Aut ( \X _{e} )$, $\Emb (\X )=\Emb ( \X _{e} )$ and $\P (\X )=\P ( \X _{e} )$.
Hence $\X$ is ultrahomogeneous iff $\X _{e}$ is.
\end{fac}
\dok
The proofs of (a), (b) and (c) are straightforward and we prove (d). For $x,y \in X$ we have:
$\la x,y \ra \in ((\r ^{-1})_{re})^c$ iff $\la x,y \ra \not\in \Delta _X \cup \r ^{-1}$ iff $x\neq y \land \la y,x \ra\not\in \r$
iff $x\neq y \land \neg y\r x  \land ( x\r y  \lor \neg x\r y  )$ iff
$(x\neq y \land \neg y\r x  \land  x\r y ) \lor
(x\neq y \land \neg y\r x  \land \neg x\r y  )$. Since the relation $\r$ is irreflexive  and asymmetric we have
$x\neq y \land \neg y\r x  \land  x\r y$ iff $x\r y$; thus $\la x,y \ra \in ((\r ^{-1})_{re})^c$ iff $ x\r y  \lor
(x\neq y \land \neg y\r x  \land \neg x\r y  )$ iff $\la x, y\ra \in \r _e$ and the equality $\X _e =((\X ^{-1})_{re})^c$ is proved.
Now applying (a) and (b) we obtain the remaining equalities. Let $\X$ be ultrahomogeneous and $\f \in \Pi ( \X _{e} )$. Then
$\f \in \Pi ( \X  )$ and, hence, there is $f\in \Aut (\X )$ such that $\f \subset f$ and, since $f\in \Aut (\X _e )$, we proved that the structure
$\X _e$ is ultrahomogeneous. The converse has a similar proof.
\hfill $\Box$

\begin{lem}\label{T9014}\rm
An irreflexive binary structure $\X$ is ultrahomogeneous  and complete iff
it is isomorphic to the enlargement of an  ultrahomogeneous digraph.
%$\X \cong \Y _e$, for some ultrahomogeneous digraph $\Y $.
\end{lem}
\dok
Let $\X= \la X,\r \ra$ be an irreflexive binary structure.

($\Rightarrow$) Assuming that $\X$ is ultrahomogeneous  and complete we define the binary relation $\rw$ on $X$ by
\begin{equation}\label{EQ9018}
x\rw y \Leftrightarrow x\r y \land \neg y\r x .
\end{equation}
{\it Claim 1.}
For the structure $\Y := \la X, \rw \ra$ we have:

(a) $\Pi (\X )=\Pi (\Y )$, $\Aut (\X )=\Aut (\Y )$ and $\Emb (\X )=\Emb (\Y )$;

(b) $\Y$ is an ultrahomogeneous digraph;

(c) $\P (\X )=\P (\Y )$;

(d) $\X =\Y _e$, that is, $\r =\rw _e$.

\vspace{2mm}
\noindent
{\it Proof of Claim 1.}
(a) It is sufficient to prove that for each $A\subset X$ and each injection $f:A\rightarrow X$ the following two conditions are equivalent:
\begin{equation}\label{EQ9019}
\forall x,y\in A \;\; (x\r y \Leftrightarrow f(x)\r f(y)),
\end{equation}
\begin{equation}\label{EQ9020}
\forall x,y\in A \;\; (x\rightarrow y \Leftrightarrow f(x)\rightarrow f(y)).
\end{equation}
Suppose that (\ref{EQ9019}) holds.
For $x,y\in A$, condition $x \rightarrow y$, that is $x\r y \land \neg y\r x$, is, by (\ref{EQ9019}), equivalent
to $f(x)\r f(y) \land \neg f(y)\r f(x)$, that is $f(x) \rightarrow f(y)$; so (\ref{EQ9020}) is true.

Let (\ref{EQ9020}) hold and $x,y\in A$. If $x=y$, then (\ref{EQ9019}) follows from the irreflexivity of $\r$. Otherwise, we have
$f(x)\neq f(y)$.

Now, if $\neg f(x)\r f(y)$, then, by (\ref{EQ9004}), $ f(y)\r f(x)$ and, hence,
$ f(y)\rightarrow f(x)$, which by (\ref{EQ9020}) implies $ y\rightarrow x$ and, hence, $\neg x\r y$. Thus
$x\r y \Rightarrow f(x)\r f(y)$.

If $\neg x\r y$, then by (\ref{EQ9004}) we have $ y\r x$ and, hence,
$ y\rightarrow x$, which by (\ref{EQ9020}) implies $ f(y)\rightarrow f(x )$ and, hence, $\neg f(x)\r f(y)$. Thus
$f(x)\r f(y) \Rightarrow x\r y$ and (\ref{EQ9019}) is true.

(b) If $\f \in \Pi (\Y )$, then, by (a), $\f \in \Pi (\X )$ and, since $\X$ is  ultrahomogeneous, there is $f\in \Aut (\X )$ such that
$\f \subset f$. By (a) again we have $f\in \Aut (\Y )$ and, thus, $\Y$ is an ultrahomogeneous structure.
Since the relation $\r$ is irreflexive, $\rw$ is irreflexive too and $x\rw y \land y\rw x$ would imply $x\r y$ and $\neg x\r y$; thus,
$\rw$ is an asymmetric relation and $\Y $ is a digraph.

(c) By (a), $\P (\X )=\{ f[X]:f\in \Emb (\X )\} =\{ f[X]:f\in \Emb (\Y )\}=\P (\Y)$.

(d) We prove that for each $x,y\in X$ we have $x\r y \Leftrightarrow x\rw _e y$, that is,
\begin{equation}\label{EQ9021}
x\r y \Leftrightarrow x\rw y \lor (x\neq y \land \neg x\rw y \land \neg y \rw x).
\end{equation}
Let $x\r y$. If $\neg y\r x$, then $x\rw y$ and, hence, $ x\rw _e y$. If $y\r x$, then, since $\r$ is irreflexive, $x\neq y$. Also
$ \neg x\rw y$ and $\neg y \rw x$ thus $ x\rw _e y$ again.

Let $ x\rw _e y$. If $x\rw y$, then $x\r y$ and we are done. If $\neg x\rw y$, then, by the assumption, $x\neq y$ and $\neg y\rw x$.
By (\ref{EQ9004}), $\neg x\r y$ would imply
$y\r x$ and, hence, $y\rw x$, which is not true. Thus $x\r y$ and Claim 1 is proved.
\kdok

($\Leftarrow$) W.l.o.g.\ suppose that $\Y =\la X , \rw \ra $ is an ultrahomogeneous digraph and $\X =\Y _e$ that is $\r =\rw _e$. Then for each $x,y \in X$ we have
\begin{equation}\label{EQ9022}
x\r y \Leftrightarrow x\rw y \lor (x\neq y \land \neg x\rw y \land \neg y \rw x).
\end{equation}
For a proof that $\X$ is complete we take different $x,y\in X$ and show that $x\r y$ or $y\r x$.
By (\ref{EQ9022}), if $x\rw y$ or $y\rw x$, then $x\r y$ or $y\r x$ and we are done. Otherwise we have $x\neq y \land \neg x\rw y \land \neg y \rw x$ and by (\ref{EQ9022}) again we obtain $x\r y$.

Since $\Y$ is an ultrahomogeneous digraph, by Fact \ref{T9011}(d) the structure $\X$ is ultrahomogeneous as well.
\kdok
\noindent
{\bf Proof of Theorem \ref{T9017}.}
Let $\X$ be an ultrahomogeneous structure and first suppose that $\X$ is disconnected.
 If $\X$ is irreflexive, then, by Theorem \ref{T9015}, $\X \cong \bigcup _\k \Y _e$, for some
ultrahomogeneous digraph $\Y$ and some $\k >1$. If $\X$ is reflexive, then $\X _{ir}$ is disconnected, irreflexive
and, by Fact \ref{T9011}(c), ultrahomogeneous so, by Theorem \ref{T9015}, $\X _{ir}\cong \bigcup _\k \Y _e$, which implies
$\X \cong (\bigcup _\k \Y _e)_{re}$. Now, suppose that $\X ^c$ is disconnected. By Fact  \ref{T9011}(a), $\X ^c$ is ultrahomogeneous.
If $\X ^c$ is irreflexive, by Theorem \ref{T9015} we have $\X ^c \cong \bigcup _\k \Y _e$, which implies $\X \cong (\bigcup _\k \Y _e )^c$. Finally,
If $\X ^c$ is reflexive, then $\X ^c _{ir}$ is disconnected, irreflexive
and, by Fact \ref{T9011}(c), ultrahomogeneous. So, by Theorem \ref{T9015} again, $\X ^c _{ir} \cong \bigcup _\k \Y _e$ which implies
$\X ^c \cong (\bigcup _\k \Y _e )_{re}$ and $\X \cong ((\bigcup _\k \Y _e )_{re})^c$.
\hfill $\Box$
\section{Posets of copies of ultrahomogeneous structures}\label{S5}
In this section we show that a classification of biconnected ultrahomogeneous digraphs, related to the properties of their posets of copies,
provides the corresponding classification inside a much wider class of structures.
\begin{te}\rm\label{T9019}
Let $\X$ be a reflexive or irreflexive ultrahomogeneous non biconnected binary structure
and let $\Y $ and $\k$ be the corresponding ultrahomogeneous digraph and the cardinal from Theorem \ref{T9017}.
Then

(a) $\la \P (\X ) , \subset  \ra \cong  \la \P (\Y) ,\subset \ra ^\k$ and
    $\sq \la \P (\X ) , \subset  \ra \cong (\sq \la \P (\Y) ,\subset \ra )^\k$, if $\k<\o$;

(b) $\sq \la \P (\X ) , \subset  \ra$ is  atomless, if $\k \geq \o$;

(c) $\sq \la \P (\X ) , \subset  \ra$ is $\k ^+$-closed, if $\k\geq \o$ is regular;

(d) $\sq \la \P (\X ) , \subset  \ra$ is forcing equivalent to the poset $(P(\k )/[\k ]^{<\k })^+$, if
$\k\geq \o$ is regular and $|\P (\Y )|\leq 2^\k =\k ^+$. The same holds for $\la \P (\X ) , \subset  \ra$.
\end{te}
\dok
By Theorem \ref{T9017}, the structure $\X$ is isomorphic to $\bigcup _\k \Y _e$,
$(\bigcup _\k \Y _e)_{re}$, $(\bigcup _\k \Y _e )^c$ or $((\bigcup _\k \Y _e)_{re})^c$ so, by Fact \ref{T9011}, $\P (\X )\cong \P (\bigcup _\k \Y _e)$.
Since the structure $\Y _e$ is complete it is strongly connected and the statement follows from Theorem \ref{T4024}.
The equality $\P (\Y _e)=\P (\Y )$ is proved in Fact \ref{T9011}(d).
\hfill $\Box$
\begin{te}\rm\label{T9020}
Let $\X$ be a countable reflexive or irreflexive ultrahomogeneous  binary structure. If $\X$ is not biconnected and $\Y $ and $\k $ are the corresponding
objects from Theorem \ref{T9017}, then

(i) $\P (\X )\cong \P (\Z )^n$, for some biconnected ultrahomogeneous digraph $\Z$ and some $n\geq 2$, if $\k<\o$ and $\Y$ has finitely many components;

(ii) $\sq \P (\X )$ is an atomless and $\o _1$-closed poset and, under CH, forcing equivalent to the
poset $(P(\o )/\Fin )^+$, if $\k=\o$ or $\Y$ has infinitely many components.
\end{te}
\dok
By Theorem \ref{T9017}, $\X$ is isomorphic to $\bigcup _\k \Y _e$,
$(\bigcup _\k \Y _e)_{re}$, $(\bigcup _\k \Y _e )^c$ or to $((\bigcup _\k \Y _e)_{re})^c$,
where $\Y$ is an ultrahomogeneous digraph and $2\leq \k \leq \o$.
So, by Fact \ref{T9011}, $\P (\X )\cong \P (\bigcup _\k \Y _e)$.

If $\k =\o$, then (ii) follows from (b), (c) and (d) of Theorem \ref{T9019}.

If $\k = n <\o$, then, by Theorem \ref{T4024}  and Fact \ref{T9011}(d), $\P (\X )\cong \P (\Y _e)^n \cong \P (\Y )^n$. We have two cases.

{\it Case 1}: $\Y$ is connected. Then, since $\Y$ is a digraph, $\Y ^c$ is a complete and, hence, a connected structure.
So $\Y$ is biconnected and we have (i).

{\it Case 2}: $\Y$ is disconnected. Then, if $\Y$ has finitely many components, say $\Y =\bigcup _{i<m}\Y _i$, by Lemma \ref{T9010}
the structures $\Y _i$ are isomorphic and complete and, hence strongly connected; so by Theorem \ref{T4024}(a), $\P (\Y )\cong \P (\Y _0)^m$, which implies
$\P (\X ) \cong \P (\Y )^n \cong \P (\Y _0)^{mn}$.
Since $\Y_0$ is a digraph and a complete structure  it is a tournament and, hence, a biconnected structure. So we have (i).

If $\Y$ has infinitely many components, say $\Y =\bigcup _{i<\o }\Y _i$, then, by Lemma  \ref{T9010}
the structures $\Y _i$ are isomorphic and complete and, hence, strongly connected. So by Theorem \ref{T4024}, the poset $\sq \P (\Y )$ is atomless and
$\o _1$-closed. Since $\P (\X ) \cong \P (\Y )^n $, by Fact \ref{T4042}(e) we have $\sq \P (\X )\cong  (\sq \P (\Y ))^n$ and, by Fact \ref{T4043}(a),
the poset $\sq \P (\X )$ is atomless and $\o _1$-closed. So we have (ii).
\kdok
\noindent
The countable ultrahomogeneous digraphs have been classified by Cherlin \cite{Che1,Che2}, see also \cite{Mac}. Cherlin's list includes
Schmerl's list of countable ultrahomogeneous strict partial orders \cite{Sch}:

- $\mathbb A _\omega$, a countable antichain (that is, the empty relation on $\omega$),

- $\mathbb B _n = n \times \mathbb Q$, for $n\in[ 1, \omega]$, where
$\langle i_1, q_1\rangle < \langle i_2 ,q_2\rangle \Leftrightarrow i_1=i_2 \;\land \; q_1 <_\mathbb Q q_2$,

- $\mathbb C _n=n\times \mathbb Q$, for $n\in[ 1, \omega]$, where  $\langle i_1, q_1\rangle < \langle i_2 ,q_2\rangle \Leftrightarrow  q_1 <_\mathbb Q q_2$,

- $\mathbb D$, the unique countable homogeneous universal poset (the random poset),

\noindent
and Lachlan's list of  ultrahomogeneous tournaments \cite{Lach}:

- $\Q$, the rational line,

- $\T ^\infty$, the countable universal ultrahomogeneous tournament,

- $S(2)$, the circular tournament (the local order),

\noindent
and many other digraphs. Also we recall the classification of countable ultrahomogeneous graphs given by Lachlan and Woodrow \cite{LachW}:

- $\BG _{\mu ,\nu }$, the union of $\mu$ disjoint copies of $\K _\nu$, where $\mu \nu=\o$,

- $\BG_{\Rado}$, the unique countable homogeneous universal graph, the Rado graph,

- $\H _n$, the unique countable homogeneous universal $\K _n$-free graph, for $n\geq 3$,

- the complements of these graphs.
\begin{ex}\rm \label{EX9000}
By the main result of  \cite{KurTod}, for the rational line, $\Q$, the poset of copies $\la \P (\Q ), \subset \ra$ is forcing equivalent to
the two-step iteration $\S \ast \pi$, where $\S$ is the Sacks forcing and $1_\S \Vdash `` \pi $ is a $\sigma$-closed forcing".
If the equality sh$(\S )=\aleph _1$ (implied by CH)
or PFA holds in the ground model, then in the
Sacks extension the second iterand is forcing equivalent to the poset $(P(\o )/\Fin )^+$.

The posets $\B _n$, $n\in[2,\o]$, from the Schmerl list are disconnected ultrahomogeneous digraphs (they are disjoint unions of copies of $\Q$)
and, by Theorem \ref{T9015}, the structures of the form $\bigcup _\k (\B _n)_e$ (or its other three variations given in Theorem \ref{T9015}) are ultrahomogeneous structures.
For example, by Theorem \ref{T9020} we have:

$\P (\bigcup _3 (\B _2)_e )\cong \P (\Q )^6 \equiv _{forc}(\S \ast \pi )^6$;

$\P ((\bigcup _\o (\B _2)_e )^c) $ and $\P (((\bigcup _2 (\B _\o )_e )_{re} )^c) $ are atomless $\o _1$-closed posets,
which are  forcing equivalent to the poset $(P(\o )/\Fin )^+$ under CH.
\end{ex}
\begin{ex}\rm \label{EX9001}
For a cardinal $\nu$, the empty structure of size $\nu$, $\A _\nu =\la \nu , \emptyset \ra$,
can be regarded as an (empty) digraph with $\nu$ components. Then $(\A _\nu )_e \cong \K _\nu$ and
for the graphs $\BG _{\mu ,\nu }$ from the Lachlan and Woodrow list we have $\BG _{\mu ,\nu }=\bigcup _\mu (\A _\nu )_e$.
So, for $n\in \N$, by Theorem \ref{T9020}, $\P (\BG _{\o ,n })$, $\P (\BG _{n, \o })$  and $\P (\BG _{\o , \o })$ are atomless $\o _1$-closed posets,
which are  forcing equivalent to the poset $(P(\o )/\Fin )^+$ under CH. But, by \cite{Kmaxemb} these posets are forcing equivalent to the posets
$(P(\o )/\Fin )^+$, $((P(\o )/\Fin )^+ )^n$ and $(P(\o \times \o )/(\Fin \times \Fin ) )^+$ respectively
and in some models of ZFC the last two of them  are not forcing equivalent to the poset $(P(\o )/\Fin )^+$.
For the first one see \cite{SheSpi1} and for the second see Example \ref{EX4015}.
\end{ex}
Let $\U$ denote the class of all countable reflexive or irreflexive ultrahomogeneous binary structures and let

$\CB =\{ \X \in \U : \X \mbox{ is biconnected}\}$,

$\CD =\{\X \in \U : \X \mbox{ is a digraph}\}$,

$\CD _e =\{\X _e : \X \in \CD \}$,

$\G =\{\X \in \U : \X \mbox{ is a graph}\}$,

${\mathcal T} =\{\X \in \U : \X \mbox{ is a tournament}\}$.

\begin{figure}[h]
\begin{center}
\unitlength 1mm %0.8mm % = .854pt
\linethickness{0.5pt}
\ifx\plotpoint\undefined\newsavebox{\plotpoint}\fi % GNUPLOT compatibility
%================================  Perfect =============================
\begin{picture}(110,102)(0,0)
%----------------------------- linije --------------------------
\put(5,5){\line(1,0){100}}%1
\put(55,5){\line(0,1){75}}%2
\put(15,10){\line(0,1){5}}%3
\put(15,10){\line(1,0){80}}%4
\put(95,10){\line(0,1){12}}%5
\put(5,5){\line(0,1){75}}%6
\put(105,5){\line(0,1){75}}%7
\put(30,20){\line(1,0){20}}%8
\put(30,20){\line(0,1){30}}%9
\put(50,20){\line(0,1){2}}%10
\put(15,15){\line(2,1){80}}%11
\put(50,29){\line(0,1){6}}%12
\put(95,30){\line(0,1){25}}%13
\put(50,35){\line(2,1){30}}%14
\put(100,30){\line(-2,1){80}}%15
\put(100,30){\line(0,1){23}}%16
\put(30,50){\line(1,0){5}}%17a
\put(45,50){\line(1,0){20}}%17b
\put(65,50){\line(0,1){15}}%18
\put(80,50){\line(0,1){15}}%19
\put(65,65){\line(1,0){15}}%20
\put(100,62){\line(0,1){13}}%21
\put(20,75){\line(1,0){80}}%22
\put(20,75){\line(0,-1){5}}%22'
\put(5,80){\line(1,0){20}}%23
\put(35,80){\line(1,0){38}}%24
\put(87,80){\line(1,0){18}}%25
%\put(80,50){\line(0,1){15}}%19

%----------------------------- tacke --------------------------

%\put(60,10){\circle*{1}}%1

%----------------------------- tekst -----------------------
%\small
\footnotesize
%\scriptsize
%\tiny

\put(46,14){\makebox(0,0)[cc]{${\mathbb C}_n$}}%a
\put(41,24){\makebox(0,0)[cc]{${\mathbb Q} $ }}%b
\put(50,25){\makebox(0,0)[cc]{${\mathcal T}$}}%c
%\put(85,14){\makebox(0,0)[cc]{${\mathbb B}_n$}}%d
\put(73,25){\makebox(0,0)[cc]{${\mathbb B}_n$}}%d
\put(95,25){\makebox(0,0)[cc]{${\mathcal D}$}}%e
%\put(20,37){\makebox(0,0)[cc]{${\mathbb U}$}}%f
\put(40,37){\makebox(0,0)[cc]{$({\mathbb C}_n )_e$}}%g
\put(40,50){\makebox(0,0)[cc]{${\mathcal D}_e$}}%h
\put(62,45){\makebox(0,0)[cc]{$({\mathbb B}_n )_e$}}%i
\put(88,46){\makebox(0,0)[cc]{${\mathbb A}_\o$}}%j
\put(73,53){\makebox(0,0)[cc]{${\mathbb K}_\o$}}%k
\put(100,57){\makebox(0,0)[cc]{${\mathcal G}$}}%l
\put(46,68){\makebox(0,0)[cc]{${\mathbb G}_{\mathrm{Rado}}$}}%m
\put(86,68){\makebox(0,0)[cc]{${\mathbb G}_{\o ,\o }$}}%n
\put(30,80){\makebox(0,0)[cc]{${\mathcal B}$}}%o
\put(80,80){\makebox(0,0)[cc]{${\mathcal U}\setminus {\mathcal B}$}}%p
\end{picture}
\end{center}
\vspace{-2mm}
\caption{Countable reflexive or irreflexive ultrahomogeneous binary structures}\label{F4002}
\end{figure}
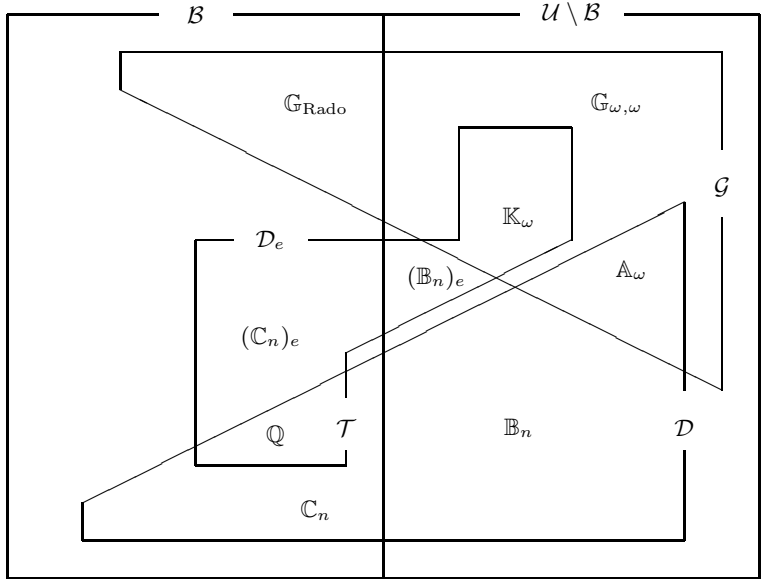

\noindent
By Lemma \ref{T9021}, the relations between these classes are displayed in Figure \ref{F4002}.
\begin{lem}\rm\label{T9021}
Let  $\Y \in \CD$. Then

(a) $\Y \in \CB$ iff $\Y$ is connected iff $\Y _e \in \CB$;

(b) $\Y  \in \CD _e$ iff $\Y$ is a tournament;

(c) $\Y \in \G$ iff $\Y =\A _\o$ iff $\Y _e =\K _\o$ iff  $\Y _e \in \G$.
\end{lem}
\dok
The first equivalence in (a) is true since $\Y ^c$ is connected, for each digraph $\Y$.
Since $\Y _e$ is connected, by Fact \ref{T9011}(d) we have $\Y _e \in \CB$ iff $(\Y _e )^c=(\Y ^{-1})_{re}$ is connected iff $\Y ^{-1}$ is connected
iff $\Y $ is connected. The statements (b) and (c) are evident.
\hfill $\Box$

\vspace{-8mm}

\noindent
By Theorem \ref{T9017} the class
$\CD $ of digraphs generates all structures from $\U \setminus \CB$ in a very simple way.
By Theorem \ref{T9020} and Fact \ref{T9011}(d), a forcing-related classification of the posets $\P (\X )$
for the structures $\X\in \CD \cap \CB$ would provide a classification for the structures $\X$ belonging to a much wider class:
$\CD \cup \CD _{re} \cup \CD _e \cup (\CD _e )_{re} \cup \U \setminus \CB$, where for a class ${\mathcal X}$ we define
${\mathcal X}_{re}=\{ \X _{re}: \X \in {\mathcal X}\}$. So, if, in addition, we obtain a corresponding classification for
$\X\in \G \cap \CB$ and hence, for $\G \cup \G _{re}$, it remains to investigate the posets $\P (\X)$ for biconnected irreflexive structures $\X$
which are not: graphs (and, hence, $\T_2\hookrightarrow \X$),
               digraphs (and, hence, $\K _2\hookrightarrow \X$),
               enlarged digraphs (and, hence, $\A _2\hookrightarrow \X$),
thus they do not have forbidden substructures of size 2.

\end{document}